\documentclass[review]{elsarticle}

\usepackage{lineno,hyperref}

\usepackage{epstopdf}
\usepackage{anysize}
\usepackage{listings}
\usepackage{graphicx}
\usepackage{float}
\usepackage{geometry}
\usepackage{booktabs}
\usepackage{array}
\usepackage{paralist}
\usepackage{verbatim}
\usepackage{subfig}
\usepackage{amsthm}
\usepackage{amssymb}
\usepackage{multirow}
\usepackage{textcomp}
\usepackage{tikz}
\usetikzlibrary{trees}
\usetikzlibrary{shapes}
\usepackage{natbib}
\usepackage{color}
\usepackage{xspace}
\usepackage{longtable}
\usepackage{lscape}
\usepackage{array}
\usepackage{amsmath}
\usepackage{float}

\usepackage{algorithm,algpseudocode}
\usepackage{caption}

\usepackage{enumitem}
\newlist{steps}{enumerate}{1}
\setlist[steps, 1]{label = Step \arabic*:}

\begin{document}

\begin{frontmatter}

\title{A Hybrid Algorithm for the Vehicle Routing Problem with AND/OR Precedence Constraints and Time Windows}

\author[mymainaddress,mysecondaryaddress]{Mina Roohnavazfar}

\author[mymainaddress]{Seyed Hamid Reza Pasandideh\corref{mycorrespondingauthor}}
\cortext[mycorrespondingauthor]{Corresponding author}
\ead{shr pasandideh@khu.ac.ir}

\author[mysecondaryaddress]{Roberto Tadei}

\address[mymainaddress]{Department of Industrial Engineering, Kharazmi University, Tehran, Iran}
\address[mysecondaryaddress]{Department of Control and Computer Engineering, Politecnico di Torino, Corso Duca degli Abruzzi, 24, Torino 10129, Italy}

\begin{abstract}
In this research, a new variant of the vehicle routing problem with time windows is addressed. The nodes associated with the customers are related to each other through AND/OR precedence constraints. The objective is minimizing the total traveling and service time. This generalization is necessary for problems where the visiting node sequence is defined according to the AND/OR relations, such as picker routing problems. We propose a Mixed Integer Linear Programming model to solve small-scale instances extended from the well-known Solomon benchmark. A meta-heuristic algorithm based on the hybridization of Iterated Local Search and Simulated Annealing approaches is developed, which can compute reasonable solutions in terms of CPU time and the accuracy of solutions. To improve the hybrid algorithm's performance, the Taguchi method is used to tune the algorithm parameters. A comprehensive computational analysis is conducted to analyze the performance of the proposed methods.
\end{abstract}


\end{frontmatter}

\section{Introduction}\label{intro}
The Vehicle Routing Problem (VRP) is a well-known problem in the Operations Research field. A given set of customers is served by a fleet of vehicles minimizing routing costs and respecting capacity constraints on vehicles. To better represent real-world problems, additional features need to be considered leading to different variants of the VRP.
A popular extension of the VRP considers time window constraints, assuming that serving a given customer must occur in a limited range. Time windows are either hard when it is prohibited to deliver outside the time interval or soft, allowing deliveries outside the boundaries against a penalty cost.

In many real industrial or social service environments where a set of tasks or services are performed, partial orders known as Precedence Constraints (PCs) are defined to represent the dependency of tasks on each other. The PCs play an essential role in a wide range of applications in various fields, such as the assembly industry, distribution of services, products, or passengers, construction projects, production scheduling, and maintenance support. In general cases, precedence constraints are represented by an arbitrary directed acyclic graph in which each arc corresponds to a precedence constraint between a pair of nodes. The precedence constraints can be classified into 'AND' and 'OR' types. The 'AND' precedence constraint means that a node can be reached only after visiting all of its predecessors. The 'OR' precedence constraint means that a node can be reached when at least one of its predecessors has been met before. The problem containing both 'AND' and 'OR' types is said to have AND/OR precedence constraints.

In the VRPs, it is also tough to suppose that the customers (nodes) are served (visited) independently of each other. However, sometimes a node must be met after visiting all its AND-type predecessors or after at least one of its OR-type predecessors.
In many real-life routing problems, serving customers may include delivery or picking up items that need to be performed in a particular order due to some reasons. They could be physical restrictions such as fragility, stackability, shape, size, and weight or preferred loading or unloading sequences to avoid extra effort on sorting the collected items at the end of the retrieving process. These restrictions may impose constraints on the visiting sequence represented by the AND/OR predecessors associated with each customer.

In this paper, we introduce a new variant of the VRP referred to as the vehicle routing problem subject to the AND/OR precedence constraints and time windows to minimize the total traveling and service time. Our interest in this problem originates from the picker routing problem, where the items are stored and retrieved manually with the order-picking sequences represented as AND/OR precedence constraints. In such a problem, to deal with various storage and replenishment policies that describe the availability of items like food or perishable products, time window limitations are defined for retrieving the items.

In the literature, no paper refers to nodes visited according to AND/OR precedence constraints to the best of our knowledge. The Pickup and Delivery problem is the most related work that includes the most straightforward form of PCs. They are represented as a single pairwise relation between pickup and delivery points. For instance, in the Dial-A-Ride problem, for any backhaul node $j$, there is a particular inhaul node $i$ where the paired PC $(i<j)$ must be met within any route. In this problem, any node has at most one predecessor that must be visited before it. However, in this work, we propose AND/OR PCs.
We formulate the problem as a Mixed Integer Linear Programming (MILP) model to optimally solve small-sized instances. A meta-heuristic algorithm based on Iterated Local Search (ILS) and Simulated Annealing (SA) approaches is developed and tuned using the Taguchi method \cite{Peace 1993} to deal with problems in a reasonable time effectively.

The contribution of this study is threefold: (i) the Vehicle Routing Problem with AND/OR Precedence Constraints and Time Windows is described; (ii) to formulate the problem, a MILP model is developed to solve small-scale instances optimally; (iii) we design and implement a meta-heuristic algorithm capable of obtaining reasonable solutions in terms of both quality and computational time.

This paper is organized as follows. Section \ref{literature} is dedicated to the literature review of the problem. The proposed problem is described and formulated in Sections \ref{problem} and \ref{model}, respectively. In Section \ref{algorithm}, we provide the meta-heuristic algorithm to address the problem. The computational results of the MILP model and the proposed algorithm are represented in Section \ref{result}. Finally, the paper ends with conclusions and interesting future research suggestions in Section \ref{conclusion}.

\section{Literature}\label{literature}

The first research on the VRP was performed by Dantzig and Ramser \cite{Dantzig and Ramser 1959}. Since then, thousands of papers on different variations of this problem have appeared in the logistics field literature.
One of the extensions of this problem, encountered in many real-life applications, is defined as the Vehicle Routing Problem with Time Windows (VRPTW). Several customers are served within predefined time windows. This problem was initially proposed by Solomon and Desorios \cite{Solomon and Desorios 1987}. The excellent surveys by Toth and Vigo \cite{Toth and Vigo 2014} and Kumar and Panneerselvam \cite{Kumar and Panneerselvam 2012} have detailed the literature of VRPTW solution approaches. Also, Dixit et al. \cite{Dixit et al. 2018} have reviewed some of the recent advancements in solving VRPTW using various meta-heuristic techniques. They are Particle Swarm Optimization (PSO), Ant Colony Optimization (ACO), Artificial Bee Colony algorithm (ABC), etc. Fan and Feng \cite{Fan and Feng 2018} applied the hill-climbing algorithm to improve the genetic algorithm to address VRPTW. They could effectively increase the speed and globally optimal solution quality of the algorithm.
Ye et al. \cite{Ye et al. 2018} addressed the time-dependent vehicle routing problem with time windows by developing a multi-type ant system algorithm hybridized with the ant colony system and the max-min ant system. The nearest neighbor selection mechanism, an insertion local search procedure, and a local optimization procedure are applied to improve the efficiency of the insertion procedure. The particle tabu search algorithm designed by Schneider et al. \cite{Schneider et al. 2018} can significantly improve the computational efficiency of the VRPTW problem and give a complete Pareto foreword. Yang et al. \cite{Yang et al. 2015} proposed the chaotic particle swarm optimization algorithm, which improves the speed, robustness, and speed of the solution of the VRPTW problem.

In many real operational research problems where a set of tasks or services are performed, it is tough to suppose that the tasks are independent of each other. In principle, the relative order between a couple of tasks is represented as a pairwise relation named Precedence Constraint (PC). These relations make the problems have a more comprehensive range of applications in various fields. The PCs arise whenever one activity or series of activities must be performed before beginning another activity or set of activities. Many examples can be mentioned, such as the assembly industry in which activities are carried out on products in different stations, logistics, construction projects, production scheduling, and maintenance support.

In the context of machine scheduling, Goldwasser and Motwani \cite{Goldwasser and Motwani 1999} derive inapproximability results for a specific single-machine scheduling problem with AND/OR precedence constraints.
Gillies and Liu \cite{Gillies and Liu 1995} addressed single and parallel machine scheduling problems to meet deadlines considering different structures of AND/OR precedence constraints. They proved NP-completeness of finding feasible schedules in many polynomially solvable settings with only AND-type precedence constraints. Moreover, they give priority-driven heuristic algorithms to minimize the completion time on a multiprocessor.
Mohring and Skutella \cite{Mohring and Skutella 2004} provided some algorithms for the more general and complex model of AND/OR precedence constraints. They showed that feasibility and questions related to generalized transitivity could be solved using essentially the same linear-time algorithm. Moreover, they discussed a natural generalization of AND/OR precedence constraints and prove that the same problems become NP-complete in this setting.
Lee et al. \cite{Lee et al. 2012} focused on flexible job-shop scheduling problems with AND/OR precedence constraints in the operations. They provided a MILP model, which can be used to compute optimal solutions for small-sized instances. They also developed a heuristic algorithm that results in a good solution for the problem regardless of its size. Moreover, a schedule builder who always gives a feasible solution and genetic and tabu search algorithms based on the proposed schedule builder were presented.
Van Den Akker et al. \cite{Van Den Akker et al. 2005} developed a solution framework which provides feasible schedules to minimize some objective function of the minimax type on a set of identical parallel machines subject to release dates, deadlines, AND/OR precedence constraints. They determined a high quality lower bound by applying column generation to the LP-relaxation.

In the context of routing problems, Moon et al. \cite{Moon et al. 2002} addressed the traveling salesman problem with precedence constraints (TSPPC). The pair-wise PCs form an order under which the nodes are visited. A genetic algorithm that involves a topological sort and a new crossover operation is proposed to solve the model. Savelsbergh and Sol \cite{Savelsbergh and Sol 1995a} presented the TSPPC model to solve the Dial-A-Ride problem where a vehicle should transport several passengers. Each passenger should be transported from a given location to a given destination. Mingozzi et al. \cite{Mingozzi et al.1997} dealt with the TSP with time windows and precedence constraints using a dynamic programming approach. Fagerholt and Christiansen \cite{Fagerholt and Christiansen 2000} considered a TSPPC with a time window to solve the bulk ship scheduling problem. The model is solved as the shortest path problem on a graph. Renaud et al. \cite{Renaud et al. 2000} proposed a heuristic model to solve the pickup and delivery TSP formulated as the TSPPC. Bredstrom and Ronnqvist \cite{Bredstrom and Ronnqvist 2008} developed a mathematical model for the combined vehicle routing and scheduling problem with time windows. The sets of pairwise synchronization and precedence constraints are considered between customer visits, independently of the vehicles. Also, they described some real-world problems to emphasize the importance of the mentioned constraints, such as homecare staff scheduling, airline scheduling, and forest operations. Bockenhauer et al. \cite{Bockenhauer et al. 2013} studied a variant of TSP in which a given subset of nodes are visited in a prescribed order in the computed Hamiltonian cycle. They presented a polynomial-time algorithm to solve the problem. Haddadene et al. \cite{Haddadene et al. 2016} modeled a home health care structure as a variant of vehicle routing problem with time windows and timing constraints. Some patients ask for more than one visit simultaneously or in given priority order. A MILP model, a greedy heuristic, two local search strategies, and three metaheuristics are proposed to solve the problem. Recently, the task assignment problem for a team of heterogeneous vehicles has been investigated in which packages are delivered to a set of dispersed customers subject to precedence constraints. Using graph theory, a lower bound on the optimal time is constructed. Integrating with a topological sorting technique, several heuristic algorithms are developed to solve it \cite{Bai et al. 2019}.

A closely related case of the VRP with PCs is the Dial-A-Ride problem, which is an exhaustively studied problem. In this problem, the pair-wise PCs are inherently represented between pickup and delivery points within a route, i.e., for any backhaul node $j$, there is a particular inhaul node $i$ where the paired PC, $(i<j)$, must be met within a route. For a comprehensive survey of the developed models, applications, and algorithms which address the Dial-A-Ride problems, the reader is referred to \cite{Ho et al. 2018}, \cite{Molenbruch et al. 2017}, and \cite{Cordeau and Laporte 2007}.

The order-picking problem is one of the main applications of PCs in the context of routing problems. However, little works in order-picking problems have focused on PCs. Zulj et al. \cite{Zulj et al. 2018} considered the PCs in a warehouse of a German manufacturer of household products, where heavy items are not allowed to be stored on top of delicate items to prevent damage to the delicate items. To avoid the sorting effort at the end of the order-picking process, they propose a picker-routing strategy respecting the precedence constraints. An exact algorithm based on dynamic programming is used to evaluate the strategy and compared with the simple s-shape routing strategy.
Dekker et al. \cite{Dekker et al. 2004} investigated combinations of storage assignment strategies and routing heuristics for a real case arising in a warehouse of a wholesaler of tools and garden equipment. A guideline has to be considered indicating that fragile products have to be picked last.
Matusiak et al. \cite{Matusiak et al. 2014} presented a simulated annealing method to address the joint order batching and precedence-constrained picker-routing problem in a warehouse with multiple depots.
The shortest path through the warehouse is determined using the exact algorithm developed by Hart et al. (1968).
Chabot et al. \cite{Chabot et al. 2017} introduced the order-picking routing problem underweight, fragility, and category constraints. They propose a capacity-index and two-index vehicle-flow formulations as well as four heuristics to solve the problem. Furthermore, a branch-and-cut algorithm is applied to solve the two mathematical models.

The background study shows no available research that can cover AND/OR precedence constraints in a capacitated vehicle routing problem. The proposed PCs are satisfied inside the route performed by each vehicle.

\section{Problem description}\label{problem}
The proposed problem is defined as a capacitated vehicle routing problem. A set of nodes are visited using a fleet of homogeneous vehicles available at the depot in time zero, with capacity $Q$.
Each vehicle can make one single trip during the planning time horizon. The problem can be represented on a directed graph $G=(N,A)$, where $N=\{0,1,...,n,\acute{0}\}$ is a set of geographically located nodes including the depot (node $0$) and a dummy depot (node $\acute{0}$), and $A=\{(i,j) | i,j \in N\}$ is a set of arcs. Each arc $(i,j) \in A$ is defined by a traveling time $t_{ij}$. For each node $i \in N \setminus\{0,\acute{0}\}$, demand $q_i$, service time $s_i$, and a (hard) time window $[e_i,l_i]$ are given where $e_i$ is the earliest possible arrival time and $l_i$ is the latest possible one. Arriving at node $i$ before $e_i$ leads to a waiting time at this node. On the other side, late arrival at the node (after $l_i$) is not allowed. A time horizon $T$ is given and establishes the working day. It can be viewed as a time window $[e_0, l_0] = [0, T]$ associated with the depot, which means the routes cannot start before $e_0$ and must be back to the depot up to time $l_0$.

Moreover, the nodes are related together by defining AND/OR precedence constraints. These relations must be met among the nodes visited by each vehicle. Given node $i$, two sets denoted by $AND_i$ and $OR_i$ are defined, including the AND-type and OR-type predecessors of node $i$, respectively. All predecessors $j \in AND_i$, visited by the same vehicle as node $i$, must be served in any positions before node $i$ on the trip. Regarding OR-type precedence constraints, at least one of the predecessors $j \in OR_i$ needs to be served before node $i$ by the same vehicle.

\begin{figure}[!ht]
\centering
\includegraphics[width=0.5\textwidth]{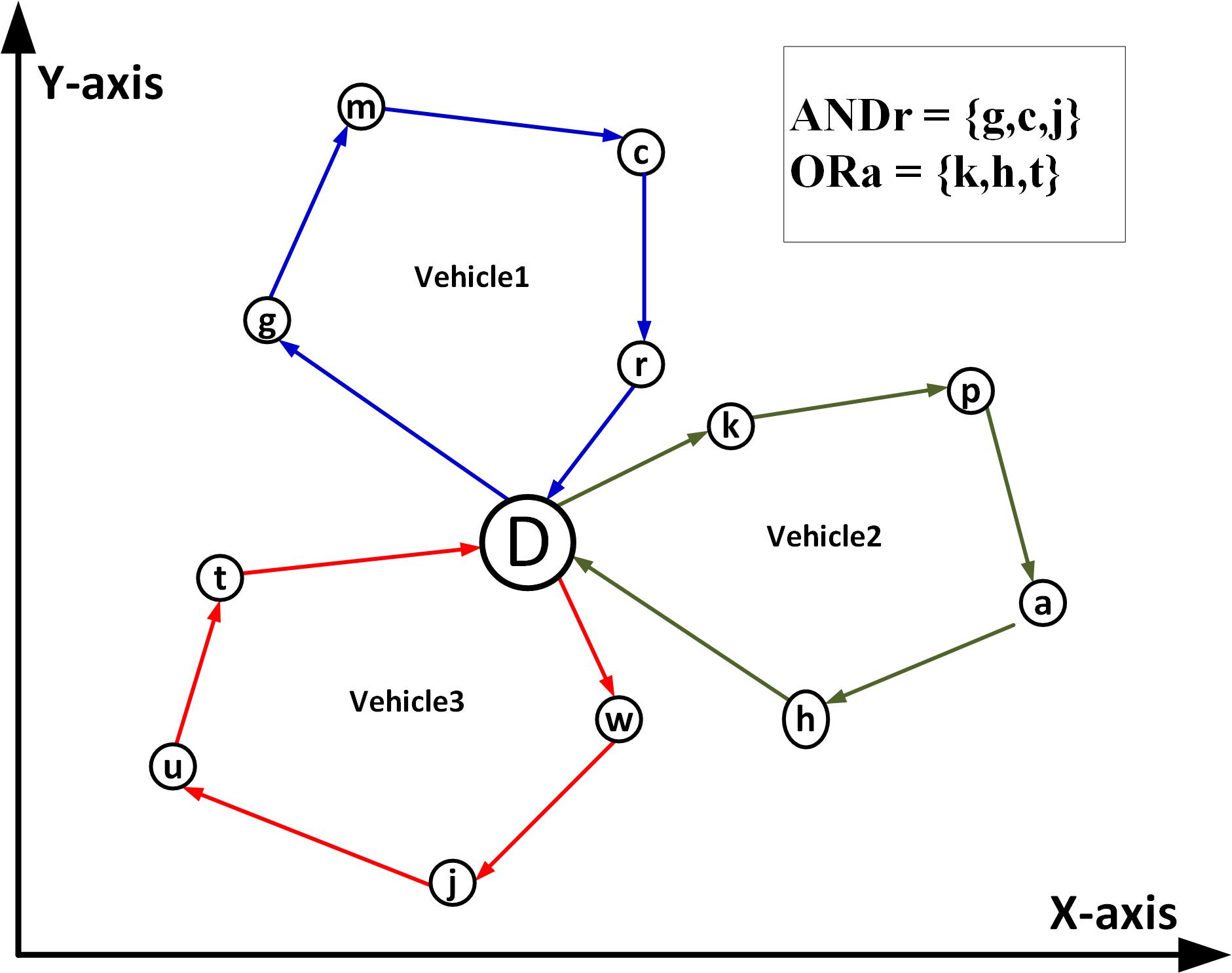}
\caption{Illustration of an example} \label{fig1}
\end{figure}

As an example, let’s consider the feasible solution depicted in Figure \ref{fig1}. It can be seen that three vehicles are used to visit all the nodes. According to the defined PCs, the set of AND-type predecessors of node $r$ includes $\{g,c,j\}$. As depicted, nodes $g$ and $c$ are visited before node $r$ using vehicle 1, while node $j$ is not, as it is not assigned to the same vehicle as node $r$.
Concerning OR-type predecessors of node $a$, node $k$ is visited before node $a$, while the corresponding precedence constraints $h < a$ and $t < a$ are not met. This is due to the definition of OR-type PCs, which denotes that given a node, at least one of its OR-type predecessors needs to be met before.

The proposed problem calls for determining an optimal trip of each vehicle to minimize the total traveling and service time to serve all the customers through network edges.

The assumptions of this problem are as follows:

\begin{itemize}
\item each vehicle starts and ends the trip at the depot;
\item each node is visited exactly once within its time window;
\item the vehicles must be back to the depot up to time T;
\item maximum number of used vehicles is $|K|$;
\item vehicle capacity $Q$ cannot be exceeded in each trip;
\item the arrival time of each node must meet the time window limitations;
\item AND/OR PCs are met between the nodes within each route.
\end{itemize}

In the following section, the proposed notations and the MILP model are presented to address the problem.

\section{Mathematical model} \label{model}
Let us introduce the following notation:
\begin{itemize}
\item $N=\{0,1,2,...,n,\acute{0}\}$: set of nodes, (each vehicle starts the route from depot (node 0) and ends to the dummy depot (node $\acute{0}$));
\item $t_{ij}$: travelling time from node $i$ to node $j$;
\item $T$: time horizon;
\item $Q$: vehicle capacity;
\item $q_i$: demand of node $i$;
\item $s_i$: service time of node $i$;
\item $K=\{1,2,...,f\}$: set of available homogenous vehicles;
\item $AND_i$: set of AND-type predecessors of node $i$;
\item $OR_i$: set of OR-type predecessors of node $i$;
\item $[e_i,l_i]$: time window associated with the arrival time of node $i$, (e.i., arriving earlier than $e_i$ introduces a waiting time at node $i$; arriving after $l_i$ leads to infeasibility);
\item $M$: an arbitrary large constant.
\end{itemize}

The variables and the proposed MILP model are presented as follows:\\

\begin{itemize}
\item $y_{ijk}$: binary variable takes value $1$ if node $i$ is visited before node $j$ (not necessarily immediately) by vehicle $k$, $0$ otherwise;
\item $z_{ik}$: binary variable takes value $1$ if node $i$ is visited by vehicle $k$, $0$ otherwise;
\item $u_{k}$: binary variable takes value $1$ if vehicle $k$ is used;
\item $a_{ik}$: continuous variable indicating the arrival time of node $i$ visited by vehicle $k$;
\item $c_{k}$: continuous variable indicating the completion time of the route performed by vehicle $k$.
\end{itemize}

The MILP model becomes:
\begin{equation}\label{obj}
\begin{aligned}
\textit{Minimize} && Z= \sum_{k=1}^{f} c_{k}
\end{aligned}
\end{equation}
subject to
\begin{equation} \label{con1}
\begin{aligned}
\sum_{k=1}^{f} z_{ik}=1 &&&&&& \forall i \in N \setminus\{0,\acute{0}\},
\end{aligned}
\end{equation}
\begin{equation} \label{con2}
\begin{aligned}
\sum_{i=1}^{n} q_i. z_{ik} \leq Q. u_{k} &&&&&& \forall k \in K,
\end{aligned}
\end{equation}
\begin{equation} \label{con3}
\begin{aligned}
z_{ik}=u_{k} &&&&&& \forall i \in \{0,\acute{0}\}, k \in K,
\end{aligned}
\end{equation}
\begin{equation} \label{con4}
\begin{aligned}
u_{k} \leq u_{\acute{k}} &&&&&& \forall k,\acute{k} \in K_{k> \acute{k}},
\end{aligned}
\end{equation}
\begin{equation}\label{con5}
\begin{aligned}
a_{0k}= 0 &&&&&& \forall k \in K,
\end{aligned}
\end{equation}
\begin{equation}\label{con6}
\begin{aligned}
a_{ik} \leq M.z_{ik} &&&&&& \forall i \in N \setminus\{0,\acute{0}\}, k \in K,
\end{aligned}
\end{equation}
\begin{equation}\label{con7}
\begin{aligned}
c_k \geq a_{\acute{0}k} &&&&&& \forall k \in K,
\end{aligned}
\end{equation}
\begin{equation}\label{con8}
\begin{aligned}
a_{ik} \geq e_i . z_{ik} &&&&&& \forall i \in N \setminus\{0\}, k \in K,
\end{aligned}
\end{equation}
\begin{equation}\label{con9}
\begin{aligned}
a_{ik} \leq l_i . z_{ik} &&&&&& \forall i \in N \setminus\{0\}, k \in K,
\end{aligned}
\end{equation}
\begin{equation}\label{con10}
\begin{aligned}
a_{jk}+ M(1-y_{ijk})\geq a_{ik}+t_{ij}+s_i-M(1-z_{ik})-M(1-z_{jk}) &&&&&& \forall i \neq j \in N, k \in K,
\end{aligned}
\end{equation}
\begin{equation}\label{con11}
\begin{aligned}
1+M(z_{ik}+z_{jk}-2) \leq y_{ijk}+y_{jik} &&&&&& \forall i \neq j \in N, k \in K,
\end{aligned}
\end{equation}
\begin{equation}\label{con12}
\begin{aligned}
z_{ik}+z_{jk} \geq 2(y_{ijk}+y_{jik}) &&&&&& \forall i \neq j \in N, k \in K,
\end{aligned}
\end{equation}
\begin{equation}\label{con13}
\begin{aligned}
y_{ijk}-1 \leq M(2- z_{ik}- z_{jk}) &&&&&& \forall i \in AND_j, k \in K,
\end{aligned}
\end{equation}
\begin{equation}\label{con14}
\begin{aligned}
1- y_{ijk} \leq M(2- z_{ik}- z_{jk}) &&&&&& \forall i \in AND_j, k \in K,
\end{aligned}
\end{equation}
\begin{equation}\label{con15}
\begin{aligned}
\sum_{i \in OR_j} y_{ijk} \geq 1+ M(z_{jk}-1) &&&&&& \forall j \in N, k \in K,
\end{aligned}
\end{equation}
\begin{equation}\label{con16}
\begin{aligned}
c_{k}, a_{ik} \geq 0 &&&&&& \forall i \in N, k \in K,
\end{aligned}
\end{equation}
\begin{equation}\label{con17}
\begin{aligned}
z_{ik}, y_{ijk}, u_k \in \{0,1\} &&&&&& \forall i,j \in N, k \in K.
\end{aligned}
\end{equation}\\
The objective function in equation (\ref{obj}) minimizes the total completion time.
Equations (\ref{con1}) guarantee that each node (except depot and dummy depot) is visited by exactly one vehicle.
The vehicle capacity is ensured for each trip by equations (\ref{con2}).
Equations (\ref{con3}), which link variables $z$ and $u$, indicate that depot and dummy depot are both assigned to all the used vehicles.
Using constraints (\ref{con4}), the sequence of used vehicles is determined so that vehicle $k$ cannot be used unless vehicle $k-1$ has been started the route before.
Equations (\ref{con5}) indicate that the depot's arrival time on each vehicle trip is equal to zero.
Constraints (\ref{con6}) relate the two assignment and arrival time variables $z$ and $a$ in a way that if node $i$ is not assigned to vehicle $k$, the corresponding arrival time of that node on that vehicle is zero.
Constraints (\ref{con7}) ensure that each vehicle trip's completion time is larger than the arrival time of the dummy depot on that vehicle.
Constraints (\ref{con8}) and (\ref{con9}) ensure that the nodes arrival time satisfy the time windows restrictions.
As represented by constraints (\ref{con10}), if nodes $i$ is visited before node $j$ on a route performed by the same vehicle $k$, the arrival time of node $j$ is larger than the sum of arrival time at node $i$, service time at this node and the travel time $t_{ij}$.
Constraints (\ref{con11}) and (\ref{con12}) ensure that if node $i$ and $j$ are both assigned to the same vehicle $k$ then node $i$ is visited either before or after node $j$ not necessarily immediately. Otherwise, the corresponding variables $y_{ijk}$ and $y_{jik}$ are equal to zero.
Constraints (\ref{con13}-\ref{con14}) ensure respecting the AND-type precedence constraints among the nodes visited by each vehicle such that if node $i$ is an AND-type predecessor of node $j$ and both nodes are assigned to the same vehicle $k$, node $i$ is visited before node $j$ which indicates $y_{ijk}=1$.
The OR-type precedence constraints are satisfied using constraints (\ref{con15}) which imply that, given node $j$ visited by vehicle $k$, at least one of its OR-type predecessors must be visited before node $j$ by the same vehicle.
Finally, constrains (\ref{con16}) and (\ref{con17}) define the continuous and binary variables, respectively.

\section{Solution approach}\label{algorithm}
In this section, our developed algorithm to address the proposed problem is introduced. Because of the high complexity of optimization problems, often exact algorithms are capable only for the smaller instances and spent a lot of computational time (see \cite{Lenstra and Rinooy Kan 1981}). In contrast, meta-heuristics can find near-optimal solutions for the instances with realistic sizes, generally with less computation time. Therefore, we concentrate on designing an effective and efficient meta-heuristic algorithm instead of exact methods.

In this research, the proposed approach is a hybridization of the Iterated Local Search (ILS) and Simulated Annealing (SA), which complements the advantages of both ILS and SA in a single optimization framework. Recently, ILS and SA have been hybridized to cope with the
search space of complex optimization problems. Experimentally, it was found that the performance of the hybrid algorithm is better than that of SA and ILS algorithms when implemented individually (see, e.g., \cite{Martin and Otto1996}, \cite{Rajalakshmi et al.2010}, and \cite{Hammouri et al.2020}). In the following, we start by describing the general schemes of the proposed algorithm. Then, the different components of this approach are represented in detail.

\subsection{General scheme}\label{General Scheme}
The general scheme of the proposed hybrid approach is represented in algorithm \ref{algorithm-General procedure}. This algorithm starts with an initial solution, denoted by $SOL$, generated from a constructive algorithm described in section \ref{Building solution}. This starting point undergoes the main loop (lines 4-22) repeated until a maximum number of iterations given by $Max_{iter}$ is reached. The loop contains four distinct parts: Perturbation Procedure (PP), Local Search (LS), check \textit{Stack} set, and updating the best solution $SOL^*$.

More specifically, the current solution is first perturbed at each iteration associated with the algorithm's destruction phase. This process starts by recognizing some target routes and modifying the solution by transferring their nodes to the not-target ones. Then, the remaining target vehicles are removed from the solution, and the corresponding nodes are gathered in a set referred to as the \textit{Stack}. More details on PP are included in Section \ref{Perturbation Procedure}.

After the destruction phase, the resulting partial solution undergoes the local search to reconstruct the solution to obtain a feasible one finally. This process is performed by randomly exploring different neighborhoods and reinserting \textit{stack} nodes to the active vehicles. In this phase, the move acceptance criterion of Simulated Annealing is incorporated with LS, which provides another way of avoiding the local optima to enhance the performance. The details on the LS are provided in Section \ref{Local Search}.

At the end of the local search loop, the current solution's feasibility is investigated by checking the \textit{stack} set that contains the nodes removed from the solution during the PP and have not been reinserted during the iterations of LS.
In case of an empty \textit{stack} set, the current solution would be feasible and can be saved as the best one if the current solution cost $F(SOL)$ is lower than the cost of the best solution $F(SOL^*)$ already found in the algorithm.
Otherwise, in the non-empty \textit{stack} set, a new vehicle may need to be added to the current partial solution. This process is implemented, taking into account parameter \textit{vehicle number}, which reports the total number of already used (active) vehicles. This parameter is updated during the algorithm by adding or eliminating each vehicle. Constructing the new vehicle route is carried out in a way that the nodes in the \textit{stack} are randomly selected, one by one, and checked to be assigned to the new vehicle till the route gets full (which means considering all feasibility aspects, no additional \textit{stack} node can be inserted to that trip).

\begin{algorithm}[H]
\caption{General procedure of the proposed hybrid algorithm}
\label{algorithm-General procedure}
\begin{algorithmic}[1]
\State Generating an initial feasible solution SOL
\State $SOL^* \leftarrow SOL$
\State $iter \leftarrow 1$
\While {$iter < Max_{iter}$}
\State Perturbation Procedure
\State Local Search
\If {\textit{stack} is not empty}  \Comment{Check \textit{Stack} Set}
\If {\textit{vehicle number} $< |K|$}
\State add a new vehicle
\State \textit{vehicle number} $++$
\State Assign \textit{stack} nodes to the new vehicle.
\State GO TO Local Search.
\Else
\State Keep the \textit{stack}
\State Go to PP
\EndIf
\EndIf
\If {$F(SOL) < F(SOL^*)$} \Comment{Update Best Solution}
\State $SOL^* \leftarrow SOL$
\EndIf
\State $iter \leftarrow iter+1$
\EndWhile
\end{algorithmic}
\end{algorithm}

The algorithm then starts the local search process again to search the neighborhoods and give the remaining \textit{stack} nodes to be inserted into the current partial solution. In case of non-empty \textit{stack} set, if parameter \textit{vehicle number} is larger than the maximum number of vehicles (K), the algorithm preserves \textit{stack}. It goes back to the Perturbation Procedure to destroy the current partial solution and reconstruct it through the LS process. The following subsections are devoted to the different parts of the proposed algorithm.

\subsection{Building initial solution}\label{Building solution}
In this section, the procedure for generating an initial solution (SOL) is introduced. This approach constructs a solution sequentially by building the routes associated with the vehicles one after another.
This process is performed by initializing an individual vehicle that travels a route during the time horizon $[0,T]$. The nodes from a sorted list of candidates are selected to be assigned to that vehicle one by one. They are located in positions one after another to form a route up to time $T$. Whenever inserting none of the unvisited nodes to the current position leads to a feasible partial solution, the depot is used to end the current trip. A new vehicle is initialized.

The partial solution is becoming complete over the stages as the number of unvisited nodes is reducing. Finally, an initial solution is obtained, which contains all the nodes that have been allocated to several vehicles (routes) that might be larger than the maximum number of available vehicles. The nodes' assignments to the vehicles are performed, considering that the selection does not violate the time horizon, vehicle capacity, time windows, and the proposed precedence constraints. However, the initial solution's feasibility in terms of the maximum number of active vehicles is achieved during the next stages of the algorithm, as described in the next sections.

In this context, we introduce parameter $AllPRE_i$, $\forall i \in N \setminus\{0,\acute{0}\}$ representing the total number of AND/OR predecessors of node $i$. Moreover, a \textit{candidate set} is constructed made of all the unvisited nodes sorted in increasing order of parameter $AllPRE$. By sorting this set, the predecessors are more likely assigned before the successors. During the nodes insertion process, the value of parameters $AllPRE$, the size and order of \textit{candidate set} are updating. When a new vehicle is added, the \textit{candidate set} is reconstructed using the unvisited vertices and their corresponding $AllPRE$. Also, this set is repeatedly updated as a node insertion is performed. Every individual node from the beginning of \textit{candidate set} is checked whether its assignment to the current vehicle leads to a feasible solution or not. Whenever a node is found whose insertion to the current position is feasible, the movement is performed.

The feasibility of each node insertion is checked in three consecutive stages. In the first stage, the partial solution's feasibility in terms of the vehicle capacity and the time horizon is investigated. In contrast, in the second and third stages, the precedence constraints and time windows are taken into account, respectively. Two parameters are defined for each vehicle, which specify the cumulative demand of already assigned nodes to that vehicle (denoted by $CD$) and the trip's completion time (denoted by $CT$).

\begin{algorithm}
  \caption{Constructing an initial solution}
  \label{algorithm-feasible solution}
  \begin{algorithmic}[1]
  \ForAll{$i \in N \setminus \{0\}$}
  \State $AllPRE_i \gets |AND_i|+|OR_i|$
  \State $\textit{VISITED} (i) \gets 0$
  \EndFor
  \State $\textrm{candidate set} \gets N \setminus \{0\}$
  \State $\textit{vehicle number} \gets 1$
  \State $\textit{current position} \gets 1$
  \While {$\textrm{candidate set} \neq \emptyset$}
\State Sort \textrm{candidate set} in increasing order of $AllPRE_i$
\State $\textit{insertion process} \gets FALSE$

\For {$i=1$ to $|\textrm{candidate set}|$}
\State Choose node $i$ \Comment{Check feasibility of inserting node $i$}
\State feasibility feedback $=$ TRUE

\State stage 1:
\If {($CD_{\textit{current position}} >$ vehicle capacity) and ($CT_{\textit{current position}} > T$)}
\State feasibility feedback $=$ FALSE
\State GO to line 11 and choose the next node
\EndIf

\State stage 2:
\If {$AND/OR$ PCs corresponding to node $i$ are not met}
\State feasibility feedback $=$ FALSE
\State GO to line 11 and choose the next node
\EndIf

\State stage 3:
\If {($e_i < A_i$) or ($A_i > l_i$)}
\State feasibility feedback $=$ FALSE
\State GO to line 11 and choose the next node
\EndIf 

\If {feasibility feedback = TRUE}
\State Insert node $i$ to current position
\State $\textit{insertion process} \gets TRUE$
\State current position $++$
\State $\textit{VISITED} (i) \gets 1$
\State Remove node $i$ from \textrm{candidate set}
\ForAll{$j \in \textrm{candidate set}$} \Comment{Update \textrm{candidate set}}
\If{$PC(i < j)= \textit{AND/OR type PC}$}
\State $AllPRE_j--$
\EndIf
        \algstore{myalg}
  \end{algorithmic}
\end{algorithm}

\clearpage

\begin{algorithm}
  \ContinuedFloat
  \caption{Constructing an initial solution (continued)}
  \begin{algorithmic}
      \algrestore{myalg}

        \If{$PC(j < i)= \textit{AND type PC}$}
\State Remove $j$ from \textrm{candidate set}
\EndIf
\EndFor 
\State Break and GO to line 8
\EndIf
\EndFor
\If {\textit{insertion process} = FALSE}
\State Add a new vehicle
\State $\textit{vehicle number}++$
\State $\textit{current position} \gets 1$
\ForAll{$i \in N \setminus \{0\}$} \Comment{reconstruct \textrm{candidate set}}
\If {$\textit{VISITED} (i)= FALSE$}
\State Add $i$ in \textrm{candidate set}
\EndIf
\EndFor 
\State Go to line 8
\EndIf
\EndWhile
  \end{algorithmic}
\end{algorithm}

For each new node insertion, the two parameters $CD$ and $CT$ are checked whether this operation leads to a partially feasible solution in terms of the vehicle capacity and time horizon. In case of a node insertion's infeasibility, the depot is placed at the end of the corresponding route, which means the vehicle ends the trip by going back to the depot.

Other feasibility aspects include precedence constraints and time windows. For each node $j$ in the \textit{candidate set}, the feasibility of corresponding precedence constraints is investigated in terms of the AND/OR PCs by considering the two following conditions:

\begin{itemize}
\item Inserting node $j$, an AND-type predecessor of any of the already assigned nodes on the route, leads to an infeasible solution.
\item Inserting node $j$ when at least one of its OR-type predecessor has not visited yet on the route leads to an infeasible solution.
\end{itemize}

Finally, the insertion of node $j$ in the \textit{candidate set} is checked by considering the time windows limitations so that node $j$ can be inserted to the current position if its arrival time meets its associated time windows.

As a node insertion is performed, the \textrm{candidate set} needs to be updated. Given node $j$ insertion to vehicle $k$, all the AND-type predecessors of node $j$, which have not been visited yet, will be removed from the \textit{candidate set}. This is due to the definition of AND-type PCs, which indicates that a vehicle cannot visit a successor before its predecessors.
Also, all successors of node $j$, which have not been assigned to any vehicles yet, are updated in the \textrm{candidate set} by subtracting constant value one from their corresponding $AllPRE$. Then, the \textrm{candidate set} is sorted again, which gives the nodes, whose predecessors have been allocated already, the chance to be assigned sooner. This procedure stops when all the nodes are assigned to the current vehicle, or no more nodes can be inserted (due to feasibility conditions). In the second case, a new vehicle is added to the current partial solution.
The above explanations on the generation of the initial feasible solution are represented as pseudo-code in algorithm \ref{algorithm-feasible solution}.

\subsection{Perturbation procedure}\label{Perturbation Procedure}
After constructing an initial solution, the algorithm starts a loop where the current solution is first perturbed with the aim of escaping from local optima. In the perturbation phase, the vehicles whose trips are not appropriate for the proportion of total traveling and service time and the number of assigned nodes are recognized and referred to as target vehicles. This process is performed using parameter $I_{k}=\frac{\textit{Traveling and service time}_{k}}{\textit{Number of nodes}_{k}}$ representing the ratio of the total traveling and service time to the number of assigned nodes associated to vehicle $k$. The vehicles whose index $I_{k}$ are higher than and equal to $I_{threshold}$ computed as equation \eqref{I_hreshold} are considered as target ones and may need to be modified.

\begin{equation}\label{I_hreshold}
I_{threshold}=\frac{\textit{Average traveling and service time}}{\textit{Average number of assigned nodes}}
\end{equation}

In equation \eqref{I_hreshold}, the average values of total traveling and service time and several nodes are taken over all active vehicles.

After recognizing the target routes, the perturbation procedure is performed in two successive steps: pre-improvement and route removal. During the first step, every individual node in target vehicles is randomly chosen and checked to transfer to one of the not-target vehicles. Whenever a not-target vehicle is found to remove the node from its current position and insert it to the new position is feasible, it is performed. The pre-improvement step aims at emptying and finally eliminating the target trips as their nodes are transferred to the not-target ones. This process continues till transferring all the target nodes are investigated, and no improvement in terms of reducing the total number of nodes in the target trips can be achieved. Any time a movement is accepted according to the feasibility test described in section \ref{Feasibility Check}, the entire solution, including the target and non-target vehicles and their assigned nodes, is updated.

\begin{algorithm}[H]
\caption{Perturbation Procedure of the proposed hybrid algorithm}
\label{algorithm PP}
\begin{algorithmic}[1]
\State Perturbation Procedure
\For {$k=1$ to $\textit{vehicle number}$}
\State Compute $I_{k}$
\If {$I_{k} > I_{threshold}$}
\State Vehicle $k$ is a target vehicle
\EndIf
\EndFor
\ForAll {target vehicles}
\ForAll {target nodes}
\State Transferring target node to a not-target route
\EndFor
\EndFor
\State Sorting the vehicles in decreasing order of $I_{k}$
\State Eliminating the first $|\textit{vehicle number}-|K||+1$ vehicles
\State Constructing \textit{stack} set using the nodes of the removed vehicles
\end{algorithmic}
\end{algorithm}

After modifying the solution in the pre-improvement step, the current solution might contain several target vehicles. In the second phase (routes removal), the algorithm destructs the solution by eliminating the first $|\textit{vehicle number}-|K||+1$ vehicles from the set of all vehicles sorted in decreasing order of $I_{k}$. The corresponding nodes of the removed vehicles are also gathered in the \textit{stack} set. After PP, the algorithm undergoes the local search iterations, described in subsection \ref{Local Search}, where the current partial solution is improved as long as the \textit{stack} nodes are reinserted to the active vehicles or the newly added ones. The above explanations on PP are represented as pseudo-code in algorithm \ref{algorithm PP}.

\subsection{Local search}\label{Local Search}
In this section, the proposed local search approach attempts to construct a feasible solution by repeatedly searching several specified neighborhoods and reinserting the \textit{stack} nodes to the current partial solution simultaneously.
This process stops when either the \textit{stack} set is empty or the maximum number of tries made without success in improving the current solution (given by $Max_{notImp}$) is reached.
At each iteration $t$ of the LS, logical parameter $Improve$ checks whether or not any improvements in terms of the objective function or reinserting nodes from the \textit{stack} set to the current partial solution is achieved. If no improvement is obtained, the parameter $iter_{notImp}$ is updated.

Most neighborhoods used in the vehicle routing problems are based on one or usually more nodes exchanges or relocations inside or across the vehicles. In this context, we apply five neighborhoods widely used in the vehicle routing problems as follows (see \cite{Petrica et al.2014}):

\begin{itemize}
\item \textbf{Transferring within a vehicle:} where a node is transferred from its current position to another position in the same vehicle.
\item \textbf{Transferring across the vehicles:} where a node is transferred from its current position to another vehicle.
\item \textbf{Exchange within a vehicle:} where two nodes that belong to the same vehicle are exchanged.
\item \textbf{Exchange across the vehicles:} where two nodes that belong to different vehicles are exchanged.
\item \textbf{Insert a vehicle:} select a node from one vehicle and create a route associated with a new vehicle with it if parameter \textit{vehicle number} is less than $K$.
\end{itemize}

During the LS procedure, the proposed five neighborhoods of the current solution are explored iteratively in a random sequence provided by a list. For example, the following sequence $[3,1,5,2,4]$ implies that neighborhood 3 is the first and neighborhood 4 is the final one to be explored at each iteration of LS.

Preliminary tests have indicated that looking for the best feasible move needs much time with no real overall improvement of the solution cost than a first improvement strategy. It should be noted that only those moves cannot lead to an infeasible solution concerning all the feasibility constraints.

To empower the LS in searching different space regions, Simulated Annealing (SA) is implemented which allows non-improving moves and avoid being trapped in a local minima. In the SA approach, a variable threshold value named \textit{SArate} is calculated as:

\begin{equation}\label{SA-Threshold}
\textit{SArate} =\exp \left(\frac{F(\acute{s})-F(s)}{Temp} \right), \;\; F(\acute{s})\geq F(s)
\end{equation}

where $F(s)$ and $F(\acute{s})$ are the objective function values of the current and the new solutions, respectively. The initial temperature $Temp$ is exponentially decreased by a fraction $\alpha$ denoted by the cooling rate, where $Temp_{t+1}=\alpha \times Temp_{t}$. In general, the SA begins with a high-temperature value, and it gradually decreases during the search. The SA performance depends on the two factors: initial temperature value $Temp$ and the cooling rate $\alpha$ calibrated in Section \ref{sec-tuning}.

For each feasible move, the origin's value and destination vehicles before and after replacement are computed. The move is accepted if it results in the lower value in terms of objective function, i.e. ($F(\acute{s})\leq F(s)$). Otherwise, the move does not improve the current solution, and it is only accepted if a random value in the range of $(0,1)$ is less than the \textit{SArate}.

\begin{algorithm}[H]
\caption{Local Search of the proposed hybrid algorithm}
\label{algorithm LS}
\begin{algorithmic}[1]
\State Local Search
\State $t \leftarrow 1$
\State $iter_{notImp} \leftarrow 1$
\While {(\textit{stack} is not empty) and ($iter_{notImp} < Max_{notImp}$)}
\State Improve$=$FALSE
\State Generate a sequence of neighborhoods randomly
\For {$n=1$ to 5}
\State Explore the neighborhood $n^{th}$
\If {F(current solution) $>$ F(new solution)}
\State Move is accepted
\State Improve$=$TRUE
\ElsIf {$\textit{random number} < \textit{SArate}$}
\State Move is accepted
\State Improve$=$TRUE
\Else
\State Move is not accepted.
\EndIf
\EndFor
\ForAll {\textit{stack} nodes}
\State Choose a node randomly
\ForAll {vehicles}
\State Choose a vehicle randomly
\ForAll {positions in the route}
\If {Inserting the node in the position is feasible}
\State Insert the node
\State Improve$=$TRUE
\State Break and go to line 18
\EndIf
\EndFor
\EndFor
\EndFor
\If {Improve=FALSE}
\State $iter_{notImp} \leftarrow iter_{notImp}+1$
\EndIf
\State $t \leftarrow t+1$
\State $Temp \leftarrow \alpha \times Temp$
\EndWhile
\end{algorithmic}
\end{algorithm}

After navigating different search space regions, adding the \textit{stack} nodes to the current partial solution is performed at each iteration of the local search cycle. The \textit{stack} nodes are randomly chosen, one by one, and checked to be inserted into the current partial solution. Given a \textit{stack} node, a vehicle is randomly selected to be allocated that node. This vehicle is chosen from all active vehicles that have not been checked before for that particular node. Chosen the vehicle, all the possible positions are assessed to be assigned the node one after another.
Whenever a vehicle and a position inside it are found under which the node insertion is feasible, the movement from the \textit{stack} to that position is implemented. This procedure aims at emptying the \textit{stack} set to construct a feasible solution finally. The above explanations on LS are represented as pseudo-code in algorithm \ref{algorithm LS}.

\subsection{Feasibility test}\label{Feasibility Check}
Since in many heuristic algorithms, millions of movements are evaluated during the search process, their feasibility must be checked as efficiently as possible. In our proposed problem, due to various side constraints such as AND/OR precedence constraints, time windows, and vehicle capacity, checking the feasibility of moves leads to prohibitive computation time. So, a critical factor of a heuristic algorithm performance is assessing the feasibility of a solution quickly. In this section, we will describe how to handle the side constraints for the considered neighborhoods efficiently.

As introduced in section \ref{Local Search}, five different moves are proposed as nodes transferring (forward or backward) or exchanging within or across the vehicles and adding new vehicles. Each move may lead to one or more operations like node removal (deletion), node insertion, forward and/or backward transfer on the corresponding vehicles. Whenever a move's feasibility is checked, a screening procedure is performed on the partial solution to return either true or false feedback. It should be noticed that when nodes transferring or exchanges are carried out within or across the routes, the feasibility conditions are only checked for the partial solution, which includes the corresponding origin and destination trips.

The feasibility test screening procedure includes three consecutive stages corresponding to checking the vehicle capacity restriction, AND/OR precedence relationships, and the time windows associated with the nodes and the depot. During the feasibility test execution, whenever an infeasibility is detected, false feedback is returned, and the procedure avoids checking other feasibility aspects in the remaining stages of the test. The proposed order has shown good performance in our preliminary experiments since the test starts from the feasibility aspects with less computational complexity to the most.

In the first stage, a partial solution's feasibility is verified concerning vehicle capacity. To do so, given a move that includes a node insertion or replacement to a position on the specific vehicle's trip, the cumulative demand of the already assigned nodes and the new one on that particular trip cannot exceed the capacity restriction.

The screening procedure checks a partial solution's feasibility in the second stage according to the precedence constraints. The cases must be considered to check the feasibility of the proposed moves concerning the AND/OR precedence constraints for each proposed neighbourhoods are represented as:

\begin{itemize}
  \item \textbf{Transferring forward within a route:} Node $i$ placed in position $p$ is transferred to a forward position $\acute{p}$ on the route.

      case 1: node $i$ should not be an AND-type predecessor of any node between location $p$ and $\acute{p}$;

      case 2: node $i$ should not be only OR-type predecessor of any node between location $p$ and $\acute{p}$.
  \item \textbf{Transferring backward within a route:} Node $i$ placed in position $p$ is transferred to a backward position $\acute{p}$ on the route.

      case 1: nodes between location $p$ and $\acute{p}$ should not be AND-type predecessor of node $i$;

      case 2: nodes between location $p$ and $\acute{p}$ should not be only OR-type predecessor of node $i$.
  \item \textbf{Transferring across routes:} Node $i$ placed in position $p$ on route $k$ is transferred to position $\acute{p}$ on route $\acute{k}$.

      case 1: node $i$ should not be only OR-type predecessor of any nodes placed after position $p$ on route $k$ ;

      case 2: node $i$ should not be AND-type predecessor of any node placed before position $\acute{p}$ on route $\acute{k}$;

      case 3: node $i$ should not be AND-type successor of any node placed after position $\acute{p}$ on route $\acute{k}$.
 \item \textbf{Exchange within a route:} Node $i$ placed in position $p$ is exchanged with node $j$ located in position $p > p$ on the same route.

     case 1: node $i$ should not be the AND-type predecessor of node $j$.

     case 2: node $i$ should not be AND-type predecessor of any node placed between position $p$ and $\acute{p}$ on the route;

     case 3: node $i$ should not be only OR-type predecessor of any node placed between position $p$ and $\acute{p}$ on the route;

     case 4: nodes between location $p$ and $\acute{p}$ should not be AND-type predecessor of node $j$;

     case 5: nodes between location $p$ and $\acute{p}$ should not be only OR-type predecessor of node $j$.
  \item \textbf{Exchange across routes:} Node $i$ placed in position $p$ on route $k$ is exchanged with node $j$ located in position $\acute{p}$ on route $\acute{k}$.

      case 1: node $i$ should not be AND-type predecessor of any node placed before position $\acute{p}$ on route $\acute{k}$;

      case 2: node $i$ should not be AND-type successor of any node placed after position $\acute{p}$ on route $\acute{k}$;

      case 3: node $i$ should not be only OR-type predecessor of any node placed after position $p$ on route $k$;

      case 4: node $j$ should not be AND-type predecessor of any node placed before position $p$ on route $k$;

      case 5: node $j$ should not be AND-type successor of any node placed after position $p$ on route $k$;

      case 6: node $j$ should not be only OR-type predecessor of any node placed after position $\acute{p}$ on route $\acute{k}$.
  \item \textbf{Insert a vehicle:} Node $i$ placed in position $p$ on route $k$ is removed and create a new route $\acute{k}$.

      case 1: node $i$ should not be only OR-type predecessor of any node placed after location $p$ on route $k$;

      case 2: node $i$ should not have any OR-type predecessor.
\end{itemize}

In the third stage, the partial solution's feasibility is checked concerning the time windows limitations. This procedure is based on the principle of push backward and forward proposed by Kindervater and Savelsberg \cite{Kindervater and Savelsberg 1997} for the Vehicle Routing Problems with Time Windows (VRPTW) where the moves are direction preserving.

Consider a path $(u,u+1,...,v)$ with associated arrival times of the nodes. Lets suppose that the arrival time of the first node in the path is decreased. This defines a push backward as:

\begin{equation}\label{push backward1}
B_u=A_u-A_u^{new},
\end{equation}

Where $A_u$ and $A_u^{new}$ define the current and new arrival time at vertex $u$. The push backward at the next vertex on the path is calculated as:

\begin{equation}\label{push backward2}
B_{u+1}=\min \{B_u, A_{u+1}-e_{u+1}\}.
\end{equation}

As long as $B_k > 0$, all nodes on the path remain feasible, and their associated arrival times need to be adjusted sequentially for $k=u,...,v$.

\begin{algorithm}[H]
\caption{Feasibility test in terms of Time Windows}
\label{test TW}
\begin{algorithmic}
\State \textbf{Transferring forward within a route:} Node $i$ placed in position $p$ is transferred to a forward position $\acute{p}$ on the route. The direct successor of node $i$ is denoted by $\sigma(i)$.
\State step 1: compute the values of push backward for the nodes placed between positions $p+1$ and
$\acute{p}-1$ and update their associated arrival times;
\State step 2: compute arrival time of node $i$ at new position $\acute{p}$;
\State step 3: compute arrival time of node $\sigma(i)$ and determine the type of push for the successors of node $i$;
\State step 4: compute the arrival times of the successors of node $i$ placed after position $\acute{p}+1$;
\State step 5: check if the new arrival times meet the time windows limitations;
\State step 6: check if the depot time window is met.
\item[]
\State \textbf{Transferring backward within a route:} Node $i$ placed in position $p$ is transferred to a backward position $\acute{p}$ on the route. The direct successor of node $i$ is denoted by $\sigma(i)$.
\State step 1: compute arrival time of node $i$ at new position $\acute{p}$;
\State step 2: compute the values of push forwards for the nodes placed between positions $\acute{p}+1$ and $p-1$ and update their associated arrival times;
\State step 3: compute arrival time of node $\sigma(i)$ and see the push is forward or backward. Then, update arrival times for the successors of node $i$;
\State step 4: check if the new arrival times meet the time windows limitations;
\State step 5: check if the depot time window is met.
\item[]
\State \textbf{Transferring across routes:} Node $i$ placed in position $p$ on route $k$ is transferred to position $\acute{p}$ on route $\acute{k}$.
\State step 1: update the arrival times of the nodes placed after position $p+1$ on route $k$ using the backward push values;
\State step 2: compute arrival time of node i in position $\acute{p}$ on route $\acute{k}$;
\State step 3: update the arrival times of the nodes placed after position $\acute{p}+1$ on route $\acute{k}$ using the push forward values;
\State step 4: check if the new arrival times meet the time windows limitations;
\State step 5: check if the depot time window is met.
\item[]
\State \textbf{Exchange within a route:} Node $i$ placed in position $p$ is exchanged with node $j$ located in position $p > p$ on the same route. The direct successor of node $i$ and $j$ on route $k$ are denoted by $\sigma(i)$ and $\sigma(j)$, respectively.
\State step 1: compute arrival time of node $j$ in new position $p$;
\State step 2: compute arrival time of node $\sigma(i)$ and determine the type of push;
\State step 3: compute the arrival times of the nodes placed between positions $p+2$ and $\acute{p}-1$ using the push values;
\State step 4: compute arrival time of node $i$ in new position $\acute{p}$;
\State step 5: compute arrival time of node $\sigma(j)$ and determine the type of push;
\State step 6: compute the arrival times of the nodes placed after positions $\acute{p}+2$ using the push values;
\State step 7: check if the new arrival times meet the time windows limitations;
\State step 8: check if the depot time window is met.

        \algstore{myalg}
  \end{algorithmic}
\end{algorithm}

\clearpage

\begin{algorithm}
  \ContinuedFloat
  \caption{Feasibility test in terms of Time Windows (continued)}\label{test TW}
  \begin{algorithmic}
      \algrestore{myalg}

\State \textbf{Exchange across routes:} Node $i$ placed in position $p$ on route $k$ is exchanged with node $j$ located in position $\acute{p}$ on route $\acute{k}$. The direct successor of node $i$ and $j$ on route $k$ are denoted by $\sigma(i)$ and $\sigma(j)$, respectively.
\State step 1: compute arrival time of node $j$ in position $p$ on route $k$;
\State step 2: compute arrival time of node $\sigma(i)$ on vehicle $k$ and determine the type of push for the next nodes on that route;
\State step 3: compute the arrival times of the nodes placed after positions $p+1$ using the push values;
\State step 4: compute arrival time of node $i$ in position $\acute{p}$ on route $\acute{k}$;
\State step 5: compute arrival time of node $\sigma(j)$ on the vehicle $\acute{k}$ and determine the type of push for the next nodes on that vehicle;
\State step 6: compute the arrival times of the nodes placed after positions $\acute{p}+1$ using the push values;
\State step 7: check if the new arrival times meet the time windows limitations;
\State step 8: check if the depot time window is met.
\item[]
\State \textbf{Insert a vehicle:} Node $i$ placed in position $p$ on route $k$ is removed and create a new route $\acute{k}$.
\State step 1: update the arrival times of the nodes placed after position $p+1$ on route $k$ using the backward push values;
\State step 2: compute arrival time of node $i$ in the first position on the route performed by vehicle $\acute{k}$;
\State step 3: check if the new arrival times meet the time windows limitations;
\State step 4: check if the depot time window is met.
  \end{algorithmic}
\end{algorithm}

Similarly, when the arrival time of the first node on path $(u,u+1,...,v)$ is postponed, a push forward is defined

\begin{equation}\label{push forward1}
F_u=A_u^{new}-A_u.
\end{equation}

The push forward at the next node on the path is calculated by

\begin{equation}\label{push forward2}
F_{u+1}=\min \{F_u-W_{u+1},0\}.
\end{equation}

where $W_{u+1}$ represents the waiting time at node $u+1$. The nodes on the path have to be checked sequentially, in such a way that if $A_k+F_k > l_k, u \leq k \leq v$, the path is no longer feasible. In case that $F_k= 0$, the path from node $k$ to node $v$ remains unchanged.

In this context, checking the feasibility of a partial solution concerning the time windows is implemented using the nodes push forward and backward for each proposed move. As mentioned before, the proposed moves may include various operations like node insertion, removal, and/or exchange with other nodes simultaneously. For example, exchanging nodes across the vehicles contains both operations of removal and insertion for two nodes on two different vehicles. So, the type of push (forward or backward) resulting from a move may not be known first. Given a movement of node $i$ on the first position of path $(u,...,v)$, it suffices only the arrival time of the direct successor of node $i$ is computed and compared with the old value. If $A_u^{new} < A_u$, the push is backward and dented by $B_u$, while in case of $A_u^{new} > A_u$ the push is forward represented by $F_u$. Then, the next nodes push on the path is sequentially computed according to \eqref{push backward2} or \eqref{push forward2} associated with the backward and forward push, respectively. Finally, the nodes' updated arrival times on the path are checked whether the associated time windows are met.

Algorithm \ref{test TW} provides the procedures implemented to check the feasibility of the proposed moves concerning the time windows associated with the related nodes and the depot (time horizon $[0,T]$).

\section{Computational results}\label{result}
In this section, we present the results of the computational experiments carried out to evaluate the MILP model and the proposed hybrid algorithm on a set of instances.
GAMS solve the MILP model. The algorithm is implemented in $C++$ on an Intel(R) Core(TM)Processor $i5-6200U$ (CPU2.30GHz) with 16 GB RAM.

In section \ref{sec-instances}, we describe the instances generated as a testbed for our assessment. In section \ref{sec-tuning}, the value of parameters involved in the algorithm is determined using the Taguchi tuning procedure. Our computational experiment results are described and commented on in Section \ref{sec results}.

\subsection{Design of experiments}\label{sec-instances}
In this section, we describe how we generate test instances for the proposed problem. Due to the novelty of the problem, no instances are available in the literature. So, we modify the well-known Solomon’s benchmark instances (see \cite{Solomon 1987}) used by most papers on vehicle routing problems with time windows.
These instances are divided into six classes obtained by combinations of two criteria. The first criterion concerns the spatial position of nodes, which includes three different options: randomly generated by Uniform distribution (denoted by R), clustered (denoted by C), and semi-clustered (denoted by RC). The second criterion is the tightness of the planning horizon, which contains two types: a short time horizon (type 1) and a long time horizon (type 2). All possible combination are therefore: "R1", "C1", "RC1", "R2", "C2" and "RC2", making a total of 56 benchmark instances. Instances are encoded as follows: C201-50 corresponds to the first instances of the class "C2", where only the first 50 customers are considered. In this work, instances with a tight planning time horizon (type 1) are discarded since the short horizon does not define a significant number of AND/OR PCs for each node.
Results are thus reported for "R2" (11 instances), "C2" (8 instances), and "RC2" (8 instances), for a total of 27 instances.

We classify the test problems into two categories in our experiments, referred to the small and large-sized instances. Both instance sets are solved using the MILP model and the proposed algorithm. The small-sized instances take the first 10, 20, and 30 nodes, while the large-sized ones take the first 40 and 50 nodes from each original instance. The nodes' locations, demands, time windows, service times, and the time horizon are set as in the original Solomon’s instances. Similar to the previous literature, the travel time is the same as the Euclidean distance between two node locations. The vehicle capacity and the maximum number of available vehicles are set to specific values, which are all empirically determined as listed in Table \ref{table value}.

\begin{table}[!ht]
\centering
\caption{Parameters setting}
\label{table value}
\small
\begin{tabular}{c c c}
\hline\noalign{\smallskip}
\textbf{Number of Nodes} & \textbf{Vehicle Capacity} & \textbf{Maximum Vehicle} \\
\hline\noalign{\smallskip}
10 & 100 & 3 \\
20 & 200 & 4 \\
30 & 200 & 4 \\
40 & 300 & 5 \\
50 & 300 & 5 \\[1ex]
\noalign{\smallskip}\hline
\end{tabular}
\end{table}

An upper triangular matrix without the diagonal called Precedence Matrix (PM) is developed to represent the precedence constraints. Each element of PM denotes whether or not a precedence relation exists between the two corresponding nodes. If node $s$ have AND-type predecessors $\{i,j,k\}$ and OR-type predecessors $\{m,n\}$, then $PM_{is}=AND$, $PM_{js}=AND$, $PM_{ks}=AND$, $PM_{ms}=OR$, and $PM_{ns}=OR$. In case that there is no PC between the two nodes, the corresponding element of the matrix is zero. The representation of precedence constraints as an upper triangular matrix leads to these relations' feasibility. It is never possible to have nodes with a smaller number than one of its predecessors.

In our proposed problem, to construct a precedence matrix, time window limitations need to be taken into account in a way that if node $i$ is a predecessor (no matter the type of PC) of node $j$, the late time of node $j$ needs to be larger than that of node $i$. Otherwise, the PC $(i<j)$ cannot be feasible in time windows.
To define a feasible precedence matrix corresponding to the time windows constraints, the set of customer nodes needs to be sorted in increasing order of the latest arrival time. The $n^{th}$ row of the matrix represents all the node's possible successors associated with that row. The $m^{th}$ column includes all the possible predecessors of the node corresponding to that column. For each column of the matrix, starting from the second column (we do not consider any predecessors for the first node) to the final one, the associated node's predecessors are generated.

To do so, we somehow adopt the scheme proposed by Derriesel and Monch \cite{Derriesel and Monch 2011} who addressed the parallel machines with sequence-dependent setup times, precedence constraints, and ready times. The precedence relations are inserted using the factor $\tau =\{0.4,0.8\}$ to evaluate PCs' impact by considering two different sizes. Given a column, if a chosen random number from $U[0,1]$ is higher than $\tau$, we do not consider any predecessors for the node associated with that column. Otherwise, the number of predecessors is chosen according to $U[0, n-1^{th}]$, where $n^{th}$ is the number of that column. Then, the predecessors are randomly selected from the set of already generated nodes $\{1^{th},...,n-1^{th}\}$. To determine the AND/OR types of precedence constraint between the randomly selected node in $\{1^{th},...,n-1^{th}\}$ and the one associated with the column, we use a random number from $U[0,1]$. The corresponding PC is an AND-type one if the chosen random number is less than $0.5$. Otherwise, the PC is an OR-type relation.

As a result, the total number of instances is $270$ which is the combinations of the type of instance $\{R2,C2,RC2\}=27$, number of nodes $\{10,20,30,40,50\}$, and the rate of precedence constraints $\{0.4,0.8\}$.
The proposed meta-heuristic algorithm is run five times over $270$ instances, and each run is stopped after 5 minutes of computation time. The upper limit of the CPU time for solving the MILP models is set to 14400 seconds.

\subsection{Tuning}\label{sec-tuning}
Because the choice of parameters has a remarkable influence on the meta-heuristic algorithms' efficiency, the Taguchi method for designing experiments is utilized to adjust the parameters.
The motivation to apply the Taguchi method in this research is that it has been recognized as an effective approach that can simultaneously consider several factors and quickly distinguish the factors with principal impacts on final solutions by performing minimal possible experiments. For more information about the Taguchi method, the interested readers can refer to \cite{Peace 1993}.

The proposed algorithm relies on four parameters, namely the maximum number of iterations $(Max_{iter})$, the maximum number of tries in LS loop are made without success to improve the current solution $(\textit{notImpMax})$, the initial temperature $(Temp)$ and the cooling rate $(\alpha)$ that need to be tuned.
The set of representative tuning instances consists of the first two instances in any combinations of the type $\{R2,C2,RC2\}$, size $N=\{10,20,30\}$, and the PC scale $\tau =\{0.4,0.8\}$ for a total number of 36 instances.
Based on some initial screening tests, three levels are selected for each parameter, as shown in Table \ref{table level}.

\begin{table}[!ht]
\centering
\caption{Parameters (factors) and their levels}
\label{table level}
\small
\begin{tabular}{c c c c}
\hline\noalign{\smallskip}
\multirow{2}{*}{\textbf{Parameters}} & \multicolumn{3}{c}{\textbf{Level}} \\
\cline{2-4}
& \textbf{1} & \textbf{2} & \textbf{3} \\
\hline\noalign{\smallskip}
$Max_{iter}$ & 500 & 900 & 1200 \\
$notImpMax$ & 50 & 70 & 90 \\
$Temp$ & 80 & 100 & 120 \\
$\alpha$ & 0.7 & 0.9 & 0.95 \\[1ex]
\noalign{\smallskip}\hline
\end{tabular}
\end{table}

Since the number of parameters and their associated levels is equal to 4 and 3, respectively, we use the orthogonal array $L_9(3^4)$. Subsequently, the number of trials and the combination of parameter levels in each trial can be specified. To enable the comparison between the objective functions of the problems, the Relative Error $(RE)$ value is calculated for each instance:

\begin{equation}\label{gap tunning}
RE(\%)=\frac{F_{sol}-F_{opt}}{F_{opt}} \times 100.
\end{equation}

where, $F_{sol}$ is the objective function value found by the algorithm and $F_{opt}$ corresponds to the optimal value for a given instance. Since $MRE$, i.e. the mean of $REs$, is to be minimized and therefore it is of “the smaller the better” category, $S/N$ ratios are obtained using following equation:

\begin{equation}\label{SN ratios}
S/N=-10 \log \left(\frac{1}{n} \sum_{i=1}^{n} F_i^2 \right).
\end{equation}

where $n$ is the number of instances, and $F_i$ is the objective function of instance $i$ obtained by the algorithm.

\begin{figure}[!ht]
\centering
\includegraphics[width=0.5\textwidth]{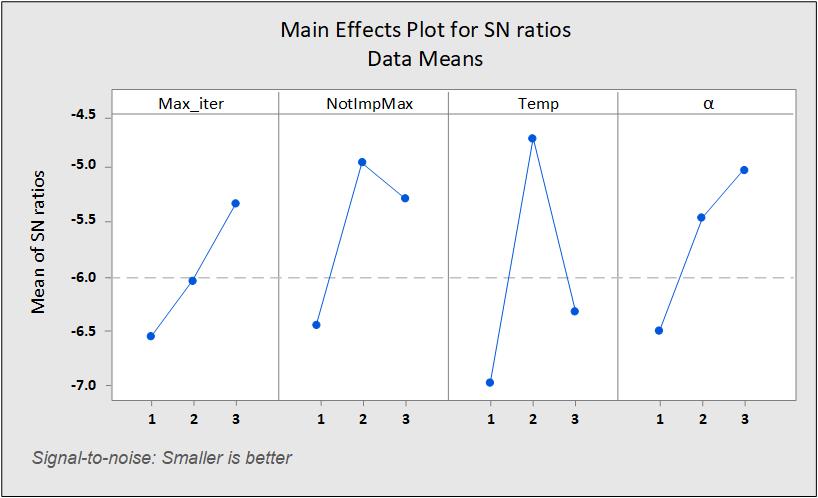}
\caption{Main effects Plots for S/N ratios} \label{fig2}
\end{figure}

\begin{figure}[!ht]
\centering
\includegraphics[width=0.5\textwidth]{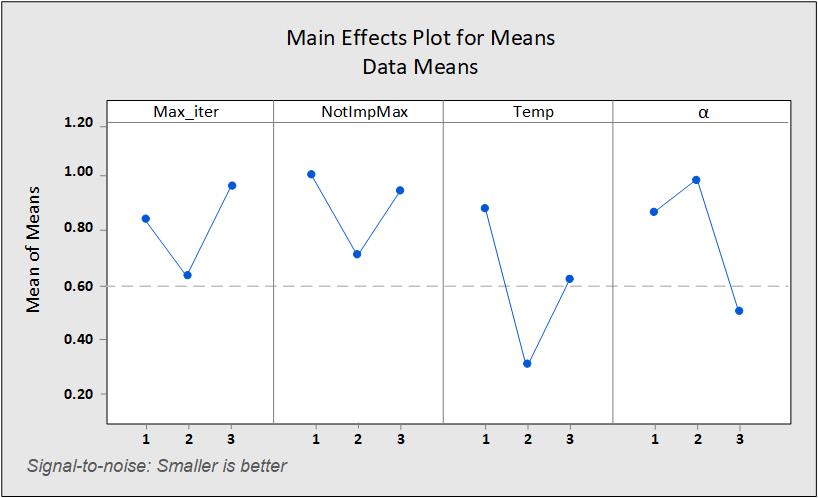}
\caption{Main effects Plots for Means} \label{fig3}
\end{figure}

Table \ref{table tuning} shows the mean value of the $S/N$ ratios and means in every level of the parameters. In this table, the parameters are ranked according to $\Delta=\max ()-\min ()$. It can be seen that the effect of $Temp$ and $Max_{iter}$ according to both signals to noise ratios and the means are the largest and the smallest, respectively, in comparison with the other factors.
Main effects plots for $S/N$ ratios and means are depicted in Figure \ref{fig2} and \ref{fig3}, respectively.
According to the signal to noise ratios, the level with maximum $S/N$ value should be selected, and then we will have $Max_{iter}=1200$, $notImpMax=70$, $Temp=100$ and $\alpha=0.95$; whereas considering the means, the level with a minimum value of mean should be chosen and therefore, we have the same results in terms of the three parameters $notImpMax$, $Temp$ and $\alpha$, with only except for $Max_{iter}$ which should be equal to $900$.
Finally, we conclude that to maintain the algorithm's robustness, the adjustment corresponding to signal to noise ratios should be selected.

\begin{table}[!ht]
\centering
\caption{Response table for the algorithm}
\label{table tuning}
\small
\begin{tabular}{c c c c c c c c c}
\hline\noalign{\smallskip}
\multirow{2}{*}{\textbf{Level}} & \multicolumn{4}{c}{\textbf{$S/N$}} & \multicolumn{4}{c}{\textbf{Mean}} \\
\cline{2-5} \cline{6-9}
&\textbf{$Max_{iter}$}&\textbf{$notImpMax$}&\textbf{$Temp$}&\textbf{$\alpha$}&
\textbf{$Max_{iter}$}&\textbf{$notImpMax$}&\textbf{$Temp$}&\textbf{$\alpha$}\\
\hline\noalign{\smallskip}
1 & -6.471 & -6.369 & -6.827 & -6.396 & 0.845 & 1.085 & 0.882 & 0.860\\
2 & -5.972 & -4.970 & -4.739 & -5.473 & 0.636 & 0.710 & 0.359 & 0.992\\
3 & -5.383 & -5.264 & -6.285 & -5.072 & 0.974 & 0.958 & 0.625 & 0.526\\
$\Delta$ & 1.088 & 1.399 & 2.088 & 1.324 & 0.338 & 0.375 & 0.523 & 0.466\\
\hline
Rank & 4 & 2 & 1 & 3 & 4 & 3 & 1 & 2\\
\noalign{\smallskip}\hline
\end{tabular}
\end{table}

\subsection{Evaluation of the MILP model and the hybrid algorithm}\label{sec results}
In this section, we summarize and discuss our experiments' results to evaluate the proposed MILP model and the algorithm performances on the different set of instances generated and described in Section \ref{sec-instances}.

Both approaches (the MILP model and the algorithm) can optimally solve all the small-sized instances. The developed algorithm's performance in terms of the computational time compared to the exact model in dealing with the first category is reported in Table \ref{table result small}. The table shows the average time (in seconds) needed to optimally solve instances for each combination of type $\{R2, RC2, C2\}$, size $N=\{10,20,30\}$, and the PC scale $\tau =\{0.4,0.8\}$ for both solution approaches.

\begin{table}[!ht]
\centering
\caption{Computational average time of the MILP model and the algorithm on small-sized instances.}
\label{table result small}
\scalebox{0.85}{
\begin{tabular}{c c c c c}
\hline\noalign{\smallskip}
\textbf{Instance} & \multicolumn{2}{c}{\textbf{MILP}} &\multicolumn{2}{c}{\textbf{Hybrid Algorithm}} \\
\cline{3-5}
\textbf{Type-N} & \textbf{$\tau =0.4$} & \textbf{$\tau =0.8$} & \textbf{$\tau =0.4$} & \textbf{$\tau =0.8$} \\
\hline\hline\noalign{\smallskip}
R2-10 & 22.51 & 16.77 & 18.92 & 24.65 \\
RC2-10 & 17.47 & 15.92 & 17.38&19.59 \\
C2-10 & 57.26 & 17.59 & 26.73& 29.13 \\
\hline
\textbf{Average:}& 25.75 & 16.76 & 21.01 & 24.46\\
\hline
R2-20 & 218.39 & 170.62 & 38.72 & 41.68 \\
RC2-20 & 378.62 & 326.35 & 27.94 & 44.80 \\
C2-20 & 583.55 &467.37 & 36.37 & 52.52 \\
\hline
\textbf{Average:}& 393.52 & 321.45 & 34.34 & 46.33\\
\hline
R2-30 & 692.07 & 583.55 & 84.99 & 127.60 \\
RC2-30 & 1073.51& 828.30 & 76.12 & 253.49 \\
C2-30 & 3802.64 &2480.73& 92.86 & 180.24 \\
\hline
\textbf{Average:}& 1856.07 & 1297.53 & 84.66 & 187.11\\
\noalign{\smallskip}\hline
\end{tabular}}
\end{table}

The first thing that should be noticed is that the proposed algorithm can find optimal solutions in less average time for all combinations except for the smallest size ($N=10$). It can be seen that the difference in time for such instances is not much as the average time over different types of size $N=10$ and PC factor $\tau=0.4$ associated to the MILP model and the algorithm are equal to $16.76$ and $21.01$, respectively.

Considering different sizes of instances, as expected, the two approaches spend more computational time to obtain optimal solutions as the number of nodes increases for all cases.

Comparing different types indicates that instances of types $R2$ and $C2$ need, respectively, the least and the most average running time to reach an optimal solution using the exact model. Instead, the developed algorithm spends average time with a small discrepancy for such instances. It means that different types of instances do not significantly affect the algorithm's computational time for instances of each size and PC factor.

Other trends can be noticed by looking at the PC scale. It seems that the larger PC scale ($\tau=0.8$) results in less complexity of the MILP model as it spends a lower average running time. On the contrary, as expected, the proposed algorithm needs more time to deal with instances with higher PC factors than those of the lower one.\\

We now want to assess the quality of solutions obtained by the second category's two approaches, including large-scale instances. The evaluation is performed by comparing their solutions before exceeding the CPU time limit of 300 and 14400 seconds associated with the algorithm and the MILP model.
The performance, in terms of percentage gap, is evaluated through the calculation of the Relative
Percentage Error (RPE) as follows:

\begin{equation}\label{RPE}
RPE(\%)=\frac{Alg_{sol}-M_{sol}}{M_{sol}} \times 100.
\end{equation}

where, $M_{sol}$ and $Alg_{sol}$ represent the feasible solution found by the MILP model and the solution obtained by the algorithm, respectively, for a given instance. It should be noticed that an optimal solution of the large-sized instances cannot be found by solving the exact model. So, the reported solution obtained by the MILP model is an upper bound with the determined optimality gap for each instance.

Tables \ref{table result2} and \ref{table result3} report the feasible solution of the MILP model and its optimality gap (denoted by OptGap(\%)), the solution found by the algorithm, as well as the factor RPE, to evaluate the performance of the algorithm concerning the feasible solution of the MILP model within the proposed time limit for instances with $N=40$ and $N=50$, respectively.

\begin{table}[!ht]
\centering
\caption{Performance results of the MILP model and the algorithm on the instances with $N=40$.}
\label{table result2}
\scalebox{0.75}{
\begin{tabular}{c| c| c c c c|c c c c}
\hline\noalign{\smallskip}
\multirow{3}{*}{\textbf{Type-N}} & \multirow{3}{*}{\textbf{Instance}} & \multicolumn{4}{c}{\textbf{$\tau =0.4$}} &\multicolumn{4}{c}{\textbf{$\tau =0.8$}} \\
\cline{3-10}
& & \multicolumn{2}{c}{\textbf{MILP}} & \textbf{Algorithm} & \multirow{2}{*}{\textbf{RPE(\%)}} & \multicolumn{2}{c}{\textbf{MILP}} &\textbf{Algorithm} & \multirow{2}{*}{\textbf{RPE(\%)}}\\
\cline{3-5} \cline{7-9}
& & \textbf{OptGap(\%)} & \textbf{Sol} & \textbf{Sol} & & \textbf{OptGap(\%)} & \textbf{Sol} & \textbf{Sol} & \\
\hline\noalign{\smallskip}
\multirow{11}{*}{R2-40}
&1 & 52.14 & 1382 & 951 & 31.19 & 48.02 & 1739 & 1226 & 29.50\\
&2 & 31.05 & 1449 & 893 & 38.37 & 34.66 & 1977 & 993 & 49.77\\
&3 & 18.12 & 1650 & 909 & 44.91 & 22.48 & 2049 & 1308 & 36.16\\
&4 & 5.89 & 1734 & 786 & 54.67 & 9.25 & 1971 & 961 & 51.24\\
&5 & 64.04 & 1250 & 991 & 20.72 & 37.61 & 1550 & 1370 & 11.61\\
&6 & 56.37 & 1484 & 852 & 42.59 & 40.07 & 1786 & 1148 & 35.72\\
&7 & 26.92 & 1451 & 1008& 30.53 & 19.22 & 1879 & 1473 & 21.61\\
&8 & 29.4 & 1359 & 661 & 51.36 & 17.39 & 1538 & 862 & 43.95\\
&9 & 47.66 & 1355 & 946 & 30.18 & 29.58 & 2150 & 1244 & 42.14\\
&10& 33.37 & 1340 & 750 & 44.03 & 31.7 & 1624 & 959 & 40.95\\
&11& 69.03 & 1240 & 837 & 32.50 & 46.44 & 1733 & 1232 & 28.91\\
\hline
\multicolumn{2}{c}{\textbf{Average:}} & 39.45& 1426.73&871.27& 38.28&30.58&1817.82&1161.45& 35.60\\
\hline
\multirow{8}{*}{RC2-40}
&1& 50.62 & 1419 & 1103 & 22.27 & 40.32 & 1792 & 1052 & 41.29 \\
&2& 36.48 & 1372 & 936 & 31.78 & 34.25 & 1885 & 1163 & 38.30 \\
&3& 54.17 & 1590 & 1108 & 30.31 & 42.66 & 1942 & 1470 & 24.30 \\
&4& 28.39 & 1648 & 950 & 42.35 & 32.72 & 1766 & 986 & 44.17 \\
&5& 41.7 & 1520 & 825 & 45.72 & 37.49 & 1827 & 992 & 45.70 \\
&6& 22.91 & 1428 & 736 & 48.46 & 24.18 & 1907 & 1373 & 28.00 \\
&7& 40.57 & 1639 & 914 & 44.23 & 12.66 & 1945 & 1328 & 31.72 \\
&8& 50.02 & 1376 & 738 & 46.37 & 47.83 & 1739 & 1014 & 41.69 \\
\hline
\multicolumn{2}{c}{\textbf{Average:}} &40.61&1499.00&913.75&38.94 &34.01&1850.38&1172.25&36.90
\\
\hline
\multirow{8}{*}{C2-40}
&1& 66.38 & 1588 & 942 & 40.68 & 51.47 & 1837 & 1059 & 42.35 \\
&2& 42.5 & 1226 & 880 & 28.22 & 36.8 & 1914 & 953 & 50.21 \\
&3& 43.99 & 1740 & 1156 & 33.56 & 41.94 & 1955 & 1274 & 34.83 \\
&4& 31.74 & 1972 & 971 & 50.76 & 29.17 & 2006 & 1368 & 31.80 \\
&5& 43.16 & 1662 & 956 & 42.48 & 44.62 & 1950 & 1076 & 44.82 \\
&6& 50.03 & 1759 & 1003 & 42.98 & 48.77 & 2173 & 1279 & 41.14 \\
&7& 25.48 & 1801 & 1149 & 36.20 & 26.83 & 1988 & 1193 & 39.99 \\
&8& 62.7 & 1839 & 944 & 48.67 & 52.9 & 2107 & 1109 & 47.37 \\
\hline
\multicolumn{2}{c}{\textbf{Average:}} &45.75&1698.38&1000.13&40.44&41.56&1991.25&1163.88&41.56
\\
\noalign{\smallskip}\hline
\end{tabular}}
\end{table}

\begin{table}[!ht]
\centering
\caption{Performance results of the MILP model and the algorithm on the instances with $N=50$.}
\label{table result3}
\scalebox{0.75}{
\begin{tabular}{c |c| c c c c |c c c c}
\hline\noalign{\smallskip}
\multirow{3}{*}{\textbf{Type-N}} & \multirow{3}{*}{\textbf{Instance}} & \multicolumn{4}{c}{\textbf{$\tau =0.4$}} &\multicolumn{4}{c}{\textbf{$\tau =0.8$}} \\
\cline{3-10}
& & \multicolumn{2}{c}{\textbf{MILP}} & \textbf{Algorithm} & \multirow{2}{*}{\textbf{RPE(\%)}} & \multicolumn{2}{c}{\textbf{MILP}} &\textbf{Algorithm} & \multirow{2}{*}{\textbf{RPE(\%)}}\\
\cline{3-5} \cline{7-9}
& & \textbf{OptGap(\%)} & \textbf{Sol} & \textbf{Sol} & & \textbf{OptGap(\%)} & \textbf{Sol} & \textbf{Sol} & \\
\hline\noalign{\smallskip}
\multirow{11}{*}{R2-50}
&1& 39.62 & 1473 & 884 & 39.99 & 18.48 & 1672 & 962 & 42.46 \\
&2& 48.5 & 1362 & 1038& 23.79 & 51.72 & 1850 & 1425 & 22.97 \\
&3& 27.31 & 1846 & 1003& 45.67 & 22.37 & 1996 & 1346 & 32.57 \\
&4& 61.17 & 1930 & 869 & 54.97 & 54.8 & 2173 & 990 & 54.44 \\
&5& 42.84 & 1393 & 920 & 33.96 & 38.51 & 1630 & 1380 & 15.34 \\
&6& 50.73 & 1560 & 1027& 34.17 & 46.04 & 1782 & 1233 & 30.81 \\
&7& 32.91 & 1528 & 940 & 38.48 & 30.18 & 1866 & 1302 & 30.23 \\
&8& 57.62 & 1402 & 841 & 40.01 & 67.82 & 1751 & 939 & 46.37 \\
&9& 40.88 & 1438 & 883 & 38.60 & 39.94 & 1849 & 1427 & 22.82 \\
&10&36.67 & 1567 & 852 & 45.63 & 41.03 & 1730 & 883 & 48.96 \\
&11&25.17 & 1346 & 733 & 45.54 & 23.74 & 1694 & 1193 & 29.57 \\
\hline
\multicolumn{2}{c}{\textbf{Average:}} &42.13&1531.36&908.18&40.07&39.51&1817.55&1189.09&34.23
\\
\hline
\multirow{8}{*}{RC2-50}
&1& 37.18 & 1583 & 1063 & 32.85 & 28.66 & 1846 & 1342 & 27.30 \\
&2& 42.93 & 1649 & 1139 & 30.93 & 38.5 & 1973 & 1274 & 35.43 \\
&3& 32.66 & 1662 & 826 & 50.30 & 30.82 & 1860 & 940 & 49.46 \\
&4& 39.03 & 1879 & 955 & 49.18 & 47.2 & 2035 & 1075 & 47.17 \\
&5& 46.82 & 1703 & 809 & 52.50 & 42.71 & 2082 & 969 & 53.46 \\
&6& 50.24 & 1534 & 783 & 48.96 & 40.9 & 1755 & 895 & 49.00 \\
&7& 56.9 & 1950 & 1094 & 43.90 & 52.52 & 2120 & 1266 & 40.28 \\
&8& 52.42 & 1526 & 1138 & 25.43 & 48.35 & 1796 & 1357 & 24.44 \\
\hline
\multicolumn{2}{c}{\textbf{Average:}} &44.77&1685.75&975.88&41.75&41.21&1933.38&1139.75&40.82
\\
\hline
\multirow{8}{*}{C2-50}
&1& 68.27 & 1837 & 1149 & 37.45 & 55.27 & 2263 & 1296 & 42.73 \\
&2& 44.62 & 1659 & 937 & 43.52 & 40.12 & 1970 & 1071 & 45.63 \\
&3& 44.09 & 1930 & 1174 & 39.17 & 43.74 & 2159 & 1289 & 40.30 \\
&4& 33.15 & 1984 & 1005 & 49.34 & 31.29 & 2324 & 1162 & 50.00 \\
&5& 44.73 & 1760 & 959 & 45.51 & 42.89 & 1973 & 1293 & 34.47 \\
&6& 66.49 & 1985 & 1340 & 32.49 & 57.18 & 2288 & 1371 & 40.08 \\
&7& 30.75 & 1874 & 967 & 48.40 & 25.36 & 2230 & 1028 & 53.90 \\
&8& 66.94 & 1936 & 1050 & 45.76 & 52.15 & 2166 & 1139 & 47.41 \\
\hline
\multicolumn{2}{c}{\textbf{Average:}} &49.88&1870.63&1072.63&42.71&43.50&2171.63&1206.13&44.32
\\
\noalign{\smallskip}\hline
\end{tabular}}
\end{table}

The first thing that should be noticed is that the algorithm's solution values are way less than the ones found by solving the MILP model for all instances of any combinations of types, sizes, and PC factors. This result shows the promising performance of the developed algorithm in dealing with such instances where the exact model cannot give optimal or even near-optimal solutions.

Considering the different types of instances, it can be seen that both OptGap and RPE associated with the instances of types $R2$ and $C2$ are the smallest and the largest value for all the cases of the two considered sizes ($N=40$ and $N=50$). It means that type $R2$ instances contain the least complexity as their optimality gap is the lower values than the two other types, while type $C2$ instances are more complex since their corresponding feasible solutions have larger OptGap values. Similar behavior can be seen for the RPE, which means that, for the instances of type $R2$, the percentage gap between the feasible solution of the MILP model and the one obtained by the algorithm is lower than those of the instances of the two other types.

Considering the two different PC factors ($\tau=0.4$ and $\tau=0.8$), a fluctuating behavior can be seen for both OptGap and RPE values in a way that they can be increased or decreased by considering the smaller or larger PC factor for each instance. The problem complexity is not entirely dependent on the amount of imposed PCs, but the generated PCs' structure may significantly affect. However, the average value of OptGap over the instances of each type and size follows a decreasing trend as the smaller PC factor is applied. The same behavior can be seen for RPE's average value over the instances except for instances of type $C2$ for both sizes.

Finally, comparing the results for the two instance sizes ($N=40$ and $N=50$) shows that, as expected, the average value of solutions obtained by the model and the algorithm and the OptGap increases as the size of instances grows. The average RPE values over the instances also follow the same described behavior with little discontinues associated with the instance of type $R2$ and PC factor $\tau=0.8$.

\section{Conclusion}\label{conclusion}
In this paper, we have studied for the first time, to the best of our knowledge, a generalization of the vehicle routing problems with time windows. The AND/OR precedence constraints are defined among the customers visited by each vehicle. This generalization comes after considering the partial orders of servicing the customers due to some physical restrictions or having the preferred loading or unloading sequence and avoiding extra effort to sort the collected items at the end of the retrieving process.
To address the problem, we have formulated it as a MILP model capable of solving only small-sized instances by spending a lot of CPU time. We have also developed a meta-heuristic algorithm as the hybridization of Iterated Local Search and Simulated Annealing approaches. The computational result of the developed algorithm highlights this approach's promising performance in CPU time and its quality.
This proves that the integration between SA and ILS can balance exploration and exploitation and thus achieve reasonable optimization results.

Future works could be devoted to different variants of routing problems like stochastic dynamic models considering AND/OR precedence constraints. Considering the solution approaches, exact methods and other algorithms can be developed to be compared with our proposed meta-heuristic algorithm in dealing with larger instances.

\end{document}